\documentclass{amsart}

\usepackage{amssymb}
\usepackage{mathrsfs}
\usepackage{amscd}

\usepackage[colorlinks,linkcolor={blue},
citecolor={blue},urlcolor={red},]{hyperref}

\theoremstyle{plain}
\newtheorem{theorem}{Theorem}[section]
\theoremstyle{remark}

\newtheorem{example}[theorem]{Example}
\theoremstyle{plain}
\newtheorem{corollary}[theorem]{Corollary}
\newtheorem{lemma}[theorem]{Lemma}
\newtheorem{proposition}[theorem]{Proposition}

\numberwithin{equation}{section}

\def\N{{\mathbb N}}

\def\R{{\mathbb R}}
\def\C{{\mathbb C}}

\newcommand{\E}{{\mathbb E}}
\renewcommand{\P}{{\mathbb P}}

\renewcommand{\d}{\delta}
\newcommand{\e}{\varepsilon}
\renewcommand{\l}{\lambda}

\newcommand{\beq}{\begin{equation}}
\newcommand{\eeq}{\end{equation}}
\newcommand{\bal}{\begin{aligned}}
\newcommand{\eal}{\end{aligned}}
\newcommand{\ben}{\begin{enumerate}}
\newcommand{\een}{\end{enumerate}}
\newcommand{\bit}{\begin{itemize}}
\newcommand{\eit}{\end{itemize}}

\newcommand{\bth}{\begin{theorem}}
\renewcommand{\eth}{\end{theorem}}
\newcommand{\bpr}{\begin{proposition}}
\newcommand{\epr}{\end{proposition}}
\newcommand{\bco}{\begin{corollary}}
\newcommand{\eco}{\end{corollary}}
\newcommand{\ble}{\begin{lemma}}
\newcommand{\ele}{\end{lemma}}
\newcommand{\bpf}{\begin{proof}}
\newcommand{\epf}{\end{proof}}
\newcommand{\bex}{\begin{example}}
\newcommand{\eex}{\end{example}}
\newcommand{\bre}{\begin{example}}
\newcommand{\ere}{\end{example}}

\renewcommand{\Re}{{\rm Re}\,}

\newcommand{\D}{{\mathscr D}}
\newcommand{\calL}{{\mathscr L}}
\newcommand{\n}{\Vert}
\newcommand{\one}{{{\bf 1}}}
\newcommand{\embed}{\hookrightarrow}
\newcommand{\s}{^*}
\renewcommand{\ss}{\subseteq}
\newcommand{\lb}{\langle}
\newcommand{\rb}{\rangle}

\newcommand{\limn}{\lim_{n\to\infty}}

\newcommand{\sumn}{\sum_{n=1}^\infty}

\newcommand{\da}{\downarrow}

\begin{document}

\title[Ornstein-Uhlenbeck semigroups]
{Norm discontinuity and spectral properties of Ornstein-Uhlenbeck semigroups}

\author{J.M.A.M. van Neerven}
\address{Delft Institute of Applied Mathematics\\
Technical University of Delft\\P.O. Box 5031\\2600 GA Delft\\The Netherlands}
\email{J.vanNeerven@math.tudelft.nl}

\author{E. Priola}
\address{Dipartimento di Matematica\\Universit\`a di Torino\\Via Carlo Alberto
10\\10123 Torino\\Italy}
\email{priola@dm.unito.it}

\date\today

\thanks{The first named author gratefully acknowledges the support by
a `VIDI subsidie' in the `Vernieuwingsimpuls' programme of the
Netherlands Organization for Scientific Research (NWO) and the
Research Training Network HPRN-CT-2002-00281. The second
named author gratefully acknowledges the support by the Italian
National Project MURST ``Equazioni di Kolmogorov''.}



\begin{abstract}
Let $E$ be a real Banach space. We study the
Ornstein-Uhlenbeck semigroup $P = \{P(t)\}_{t\ge 0}$ associated
with the Ornstein-Uhlenbeck operator
$$ Lf(x) = \tfrac12{\rm Tr}\, Q D^2 f(x) + \lb Ax, Df(x)\rb, \qquad x\in E.$$
Here $Q\in\calL(E\s,E)$ is a positive symmetric operator and $A$
is the generator of a $C_0$-semigroup $S= \{S(t)\}_{t\ge 0}$ on
$E$. Under the assumption that $P$ admits an invariant measure
$\mu_\infty$ we  prove that if $S$ is eventually compact and the
spectrum of its generator is nonempty, then
$$\n P(t)-P(s)\n_{\calL( L^1(E,\mu_\infty))} = 2 \ \ \hbox{for all $t,s\ge 0$
with $t\not=s$.}
$$
This result is new even when $E = \R^n$.
We also study the
behaviour of $P$ in the space $BUC(E)$. We show that
if $A\not=0$ there exists $t_0>0$ such that
$$\n P(t)-P(s)\n_{\calL(BUC(E))} = 2 \ \ \hbox{for all $0\le t,s\le t_0$ with $t\not=s$.}
$$
Moreover, under a nondegeneracy assumption or a strong Feller
assumption, the following dichotomy holds: either
$$ \n P(t)- P(s)\n_{\calL(BUC(E))} = 2 \ \
\hbox{for all $t,s\ge 0$, \ $t\not=s$,}
$$
or $S$ is the direct sum of a nilpotent semigroup and a
finite-dimensional periodic semigroup.  Finally we investigate
the spectrum of $L$ in the spaces $L^1 (E,\mu_\infty)$ and
$BUC(E)$.
\end{abstract}

\maketitle

\section{Introduction and Preliminaries}\label{sec:intr}

In this paper we study certain properties of the
Ornstein-Uhlenbeck semigroup in spaces of continuous functions and
integrable functions. This semigroup is associated with the
stochastic linear Cauchy problem \beq\label{eq:ACP} \left\{
\bal dU(t) & = AU(t)\,dt + B\,dW_H(t), \\
      U(0) & = x.
\eal \right. \eeq Here $A$ is assumed to be the infinitesimal
generator of a $C_0$-semigroup $S=\{S(t)\}_{t\ge 0}$ on a real
Banach space $E$, $B$ is a bounded operator from a real Hilbert
space $H$ into $E$, $W_H = \{W_H(t)\}_{t\ge 0}$ is an
$H$-cylindrical Brownian motion, and $x\in E$ is an initial value.
As is well known, the above problem admits a unique weak solution
if and only if for all $t\ge 0$ there exists a centred Gaussian
Radon measure $\mu_t$ on $E$ whose covariance operator
$Q_t\in\calL(E\s,E)$ is given by
$$
\lb Q_t x\s,
y\s\rb = \int_0^t \lb S(s)BB\s S\s(s)x\s,y\s\rb\,ds,\qquad x\s,y\s\in E\s,
$$
where $E\s$ denotes the topological dual of $E$. Under this assumption
the solution $U= \{U(t,x)\}_{t\ge 0}$ of \eqref{eq:ACP} is given
by the stochastic It\^o integral
$$ U(t,x) =  S(t)x + \int_0^t S(t-s) B \,dW_H(s),$$
see \cite{BN,DZ,NW}. For more information on
Gaussian measures in infinite dimensions we
refer to  \cite{Bo,VTC}.

The {\em Ornstein-Uhlenbeck semigroup} $P=\{P(t)\}_{t\ge 0}$ associated with
$A$ and $B$ is defined on the space $C_b(E)$ of bounded real-valued continuous
functions on $E$ by
 \beq \label{v}
P(t)f(x) := \E(f(U(t,x))) = \int_E f(S(t)x+y)\,d\mu_t(y), \qquad
x\in E, \ \   f\in  C_b(E).
 \eeq
This semigroup leaves $BUC(E)$, the space of bounded real-valued
uniformly continuous functions on $E$, invariant and has been studied by
many authors \cite{CeG,DL,GK,GN,Ku,NZ,Pr}. It is well known that
$P$ fails to be strongly continuous with respect to the supremum
norm of $BUC(E)$ unless $A=0$. Therefore it is natural to
introduce the closed subspace $BUC^\circ(E)$ consisting of all
functions on which $P$ acts in a strongly continuous way. This
subspace is invariant under $P$, and the restriction $P^\circ$  of
$P$ is strongly continuous on $BUC^\circ(E)$. It is well known
that the behaviour of $P^\circ$ is quite pathological. For
instance, in the setting of a Hilbert space $E$ it was shown in
\cite{NZ} that one has
\beq\label{eq:P2} \n P^\circ(t)
-P^\circ(s)\n_{\calL(BUC^\circ(E))} = 2
 \eeq whenever
$\mu_t\perp\mu_s$, i.e., the measures $\mu_t$ and $\mu_s$ are
mutually singular. Here $\n\cdot\n_{\calL(X)}$ denotes the uniform
operator norm of the Banach space $\calL(X)$ of all bounded linear
operators on $X$. For the heat semigroup, which corresponds to the
case $A=0$, \eqref{eq:P2} was established earlier in \cite{DR}. In
Section \ref{sec:Cb} we extend this result to Banach spaces and
complement it by showing that \eqref{eq:P2} also holds whenever
$S(t)\not=S(s)$. It follows that if $A\not=0$, then there exists
$t_0>0$ such that
\beq\label{1} \n P^\circ(t)-P^\circ(s)\n_{\calL(BUC^\circ(E))}=2
\ \ \hbox{for all $0\le s,t\le t_0$, \ $t\not=s$.}
\eeq In
particular, if $A\not=0$, then $P^\circ$ always fails to be
norm continuous on $BUC^\circ(E)$ for $t>0$.
In the converse direction we show that
{for fixed $t,s\ge 0$,} \eqref{eq:P2} and $S(t)=S(s)$ imply
 $\mu_t\perp\mu_s$.
These results are used to prove the following dichotomy: either
 \beq \label{2}
\n P^\circ(t)- P^\circ(s)\n_{\calL(BUC^\circ(E))} = 2 \ \
\hbox{for all $t,s\ge 0$, \ $t\not=s$,} \eeq or $S$ is the direct
sum of a nilpotent semigroup and a finite-dimensional periodic
semigroup. Note that this result is  new even when $E $ is
finite-dimensional. The probabilistic interpretation of \eqref{1} and \eqref{2}
is that $\sup_{x \in E} \, $ $\n \mu_{t, x} - \mu_{s, x}
\n_{{\rm var}}=2$, for $t,\, s\ge 0$ with $t\not=s$, where $\mu_{t,
x}$ denotes the law of the process $U(t,x)$ which solves
\eqref{eq:ACP}, and $\n\cdot\n_{\rm var}$ is  the total variation
norm.

{\vskip 1mm} Related to the problem of norm discontinuity is the
problem of characterizing the spectrum of the generator
$L_{P^\circ}$ of $P^\circ$. For finite-dimensional spaces $E$, it
was shown in \cite{M} that if the operator $Q:= B\circ B^*$
is invertible and the spectrum
$\sigma (A)$ of $A$ is contained in $\{\l\in\C: \ \Re\l<0\}$, then
\beq\label{eq:spectrum} \sigma(L_{P^\circ}) = \C^-, \eeq where
$\C^-:= \{\l\in\C: \Re\l\le 0\}$, and every $\l\in\C$ with
$\Re\l<0$ is an eigenvalue.
 By standard results from semigroup
theory,    \eqref{eq:spectrum} already implies that   $P^\circ$
cannot be eventually norm continuous in $BUC^\circ(E)$.
Below we obtain an extension of \eqref{eq:spectrum} to the case
where $S$ is an eventually compact semigroup on a Banach space
$E$.

\medskip Let us next assume that the limit $Q_\infty :=
\lim_{t\to\infty}Q_t$ exists in the weak operator topology of
$\calL(E\s,E)$ and that there exists a centred Gaussian Radon
measure $\mu_\infty$ with covariance operator $Q_\infty$. A
sufficient condition for this is that the Gaussian Radon measures
$\mu_t$ exist and $S$ is uniformly exponentially stable; cf.
\cite[Chapter 9]{DZ}, \cite{NW2}. The measure $\mu_\infty$ is {\em
invariant} for $P$, in the sense that for all $f\in BUC(E)$ and $t\ge 0$,
 $$
\int_{E} P(t)f(x)\,d\mu_\infty(x) = \int_{E}f(x)\,d\mu_\infty(x).
$$ By a standard argument, the semigroup $P$ has a unique
extension to a strongly continuous contract\-ion semigroup on
$L^p(E,\mu_\infty)$  for all $p\in [1,\infty)$. For $p\in
(1,\infty)$, the behaviour of this semigroup is well understood.
We refer to \cite{CG,CG2,MPP,MPRS,Ne1}, where  the domain of the
generator, its spectrum, and analyticity properties are
characterized.

The behaviour of $P$ in $L^1(E,\mu_\infty)$ is much less well
understood. For finite-dimensional spaces $E$ it is shown in
\cite{MPP} that the $L^1(E,\mu_\infty)$-spectrum of its generator
$L_P$ equals $\C^-$. To the best of our knowledge, it is an open
problem whether this result extends to infinite dimensions.
Furthermore no $L^1$-analogue of \eqref{eq:P2} seems to be known.
In Section \ref{sec:L1} we will first show, for
finite-dimensional spaces $E$, that
$$
\n P(t)-P(s)\n_{\calL(L^1(E,\mu_\infty))}=2
$$
whenever $t>s\ge 0$.  Then we extend this result
to infinite dimensions in the setting of eventually compact semigroups $S$,
and, extending a result for $E=\R^d$ in \cite{MPP},
we prove that the spectrum of $L_P$ equals $\C^-$.

{\vskip 2mm} Our approach is based on a technique introduced by
Davies and Simon \cite{DS} which may be described as follows. If
$B_1$ and $B_2$ generate $C_0$-semigroups of contractions $T_1$
and $T_2$ on a Banach space $X$, then $B_1$ {\em belongs to the
limit class} of $B_2$ if there exists a sequence of invertible
isometries $V_n: X\to X$ such that
$$  R(\l,B_1)x = \limn V_n^{-1}R(\l,B_2)V_nx, \qquad x\in E.
$$
 Here $R(\l,B_k) = (\lambda - B_k)^{-1}$, $k=1,2$. This is
 equivalent to require that, for each $t>0$,
$$
 T_1 x = \limn V_n^{-1} T_2 V_n x, \qquad x\in E.
$$
In this situation one has
$$ \n T_2(t) - T_2(s)\n_{\calL(X)} \ge   \n T_1(t) - T_1(s)\n_{\calL(X)}, \qquad t,s\ge 0,
\;\;$$  and $$ \;\; \n R(\l,B_2)\n_{\calL(X)} \ge \n
R(\l,B_1)\n_{\calL(X)}, \qquad \l\in\varrho_\infty(B_1)
\cap\varrho_\infty(B_2),$$ where $ \varrho_\infty(B_k)$ denotes
the connected component of the resolvent set $\varrho(B_k)$
containing $+\infty$, $k=1,2$.  This technique is applied in the
situation where $B_2$ is a suitable realization of the generator
of $P$ and $B_1$ is a realization of the generator of the {\em
drift semigroup} $R$ associated with $A$. This semigroup is
defined on $C_b(E)$ by
\begin{equation}\label{3}
R(t)f(x) := f(S(t)x), \qquad x\in E, \ \  f\in C_b(E).
\end{equation}

 Throughout this paper,
a {\em Gaussian measure} is a centred Gaussian Radon measure.

\section{The Ornstein-Uhlenbeck semigroup in spaces of continuous functions}
\label{sec:Cb}

In this section we study various properties of the
Ornstein-Uhlenbeck semigroup $P$ and the drift semigroup $R$ in
the spaces $C_b(E)$ and $BUC(E)$. We denote by $\n \cdot \n $
the supremum norm.

As semigroups on $C_b(E)$, both $P$ and $R$ are strongly
continuous with respect to the mixed topology.
This topology is defined as the finest locally convex topology
in $C_b(E)$ which agrees on every norm bounded set with the
topology of uniform convergence on compact sets; see \cite{Wh,Wi}
for a detailed investigation of its properties. This topology is
complete and may be used to define the infinitesimal generators
$L_P$ and $L_R$ of $P$ and $R$ by taking, for $T=P$ or $R$,
$$
\begin{aligned}
\D(L_T) & := \Big\{ f\in C_b(E): \ \lim_{t\downarrow 0} \frac1t (T(t)f - f) \
\hbox{ exists }\Big\}, \\
L_Tf & := \lim_{t\downarrow 0} \frac1t(T(t)f - f), \qquad f\in
\D(L_T),
\end{aligned}
$$
where the limits are taken with respect to the mixed topology.
We have $f\in\D(L_T)$ if and only if the following two conditions hold:
\begin{enumerate}
\item[(i)] $\displaystyle\limsup_{t\downarrow 0}
\frac1t\|T(t)f -f\|<\infty$;
\item[(ii)] there exists a function $g\in C_b(E) $ such that for
all $x\in E$,
$$\lim_{t\downarrow 0}\frac1t(T(t)f (x)-f (x))=g(x).$$
\end{enumerate}
In this situation, $L_Tf = g$.

On a suitable core of smooth cylindrical functions, $L_P$ and
$L_R$ are given by
$$
\bal L_Pf(x) &= \tfrac12{\rm Tr}\,QD^2 f(x) + \lb Ax, Df(x)\rb,\\
 L_Rf(x) &= \lb Ax, Df(x)\rb,
\eal
$$
where `Tr' denotes the trace and $Q := BB\s$. We refer to \cite{GK,GN} for proofs
and more details. Alternative approaches to diffusion semigroups
in spaces of continuous functions may be found in
\cite{CeG,Ku,Pr}.

Both $P$ and $R$ leave the closed subspace $BUC(E)$ of $C_b(E)$
invariant, but even on this smaller space both semigroups  fail to
be strongly continuous with respect to the supremum norm,
unless $A=0$. It is easy to see, cf. \cite[Lemma 3.2]{DL}, that
the closed subspaces of $BUC(E)$ on which $P$ and $R$ act in a strongly
continuous way with respect to the supremum norm coincide. This
common subspace will be denoted by $BUC^\circ(E)$. Thus,
$$
\bal BUC^\circ(E)
& = \big\{f\in BUC(E): \ \lim_{t\downarrow 0} \n P(t)f-f\n=0\big\}
\\ & = \big\{f\in BUC(E): \ \lim_{t\downarrow 0}
\n R(t)f-f\n=0\big\}.
\eal$$
The restrictions of $P$ and $R$ to $BUC^\circ(E)$, denoted
by $P^\circ$ and $R^\circ$ respectively, are strongly continuous
with respect to the supremum norm. Their generators $L_{P^\circ}$
and $L_{R^\circ}$ are characterized as follows; see
\cite[Proposition 3.5]{DL} for a related result.

\begin{proposition} We have
$$
\bal
\D(L_{P^\circ}) &= \big\{f\in \D(L_P)\cap BUC^\circ(E): \ L_Pf \in  BUC^\circ(E)\big\},\\
\D(L_{R^\circ}) &= \big\{f\in \D(L_R)\cap BUC^\circ(E): \ L_Rf \in  BUC^\circ(E)\big\}.
\eal
$$
\end{proposition}
\bpf
Let $T=P$ or $R$.

The inclusion `$\ss$' is clear. To prove the inclusion
`$\supseteq$' let $f\in \D(L_T)\cap BUC^\circ(E)$ be such that $L_T
f\in BUC^\circ(E)$. Then,
$$
\bal
\ & \lim_{t\da 0}  \sup_{x\in E} \Big|\frac{1}{t}\big(T(t)f(x)-f(x)\big)-L_Tf(x)\Big|
\\ & \qquad = \lim_{t\da 0} \sup_{x\in E}\Big|\frac{1}{t}\int_0^t T(s)L_Tf(x)\,ds - L_Tf(x)\Big|
 =\lim_{t\da 0} \Big\n\frac{1}{t}\int_0^t T^\circ(s)L_Tf-L_Tf\,ds\Big\n =0,
\eal
$$
where the first identity is a consequence of the fact that $T$ is
strongly continuous with respect to the mixed topology. This
proves that $f\in \D(L_{T^\circ})$. \epf

We do not know whether $\D(L_{P^\circ})$ is always contained in
$\D(L_{R^\circ})$.

{\vskip 1mm} The following simple observation, cf. the proof of
\cite[Lemma 2.3]{NZ}, will be useful.

\begin{lemma}\label{lem:delta}
Let $T=P$ or $R$. For $f\in BUC(E)$ and $\d>0$ define
$$ f_\d(x):= \frac1\d \int_0^\d T(t)f(x)\,dt, \quad x\in E.$$
Then $f_\d\in BUC^\circ(E)$. Moreover, $\lim_{\d\da 0} f_\d = f$ in the mixed topology inherited from $C_b(E)$.
\end{lemma}

\begin{proof}
First note that $t\mapsto T(t)f(x)$ is continuous for all $x\in E$, and therefore the function $f_\d$ is well defined.
It is clear that $f_\d\in BUC(E)$ and $\n f_\d\n \le 1$.
For all $x\in E$ and $t\in (0,\d)$ we have
 $$|T(t)f_\d(x) - f_\d(x)| = \frac1\d\Big| \int^{\d+t}_t
T(s)f(x) \,ds-\int^{\d}_0 T(s)f(x)\,ds\Big|\le \frac{2t}{\delta} \n f \n.$$
Thus $\n T(t) f_\d - f_\d\n \le 2\d^{-1}t \n f \n$, which
shows that $f_\d\in BUC^\circ(E)$. The final statement is obvious.
\end{proof}

Obviously, if $S(t)=S(s)$ for certain $t,s\ge 0$, then  $R(t)=R(s)$. The following lemma describes what happens if  $S(t)\not=S(s)$.

\ble\label{lem:norm-disc-R-Cb} For all $t,s\ge 0$ such that
$S(t)\not=S(s)$ we have $\n
R^\circ(t)-R^\circ(s)\n_{\calL(BUC^\circ(E))} =2.$
\ele \bpf Fix
$t,s\ge 0$ such that $S(t)\not=S(s)$. We may assume that $t > s\ge
0$. Choose $x_0\s\in E\s$ such that $S\s(t)x_0\s\not=S\s(s)x_0\s$.
Noting that $S\s(s)x_0\s\not=0$ we pick $x_0\in  E$ such that $
\lb x_0,S\s(t)x_0\s\rb= 0$ and $\lb x_0,S\s(s)x_0\s\rb=\pi$. The
function $f(x) := \cos\lb x,x_0\s\rb$ defines an element of
$BUC(E)$ and we have
$$\n R(t)f - R(s)f\n \ge \big|R(t)f(x_0) - R(s)f(x_0)\big| = 2.
$$
Given $\e>0$ we choose $\d>0$ small enough such that
$$\big|R(t)f_\d(x_0) - R(s)f_\d(x_0)\big| = \big|(R(t)f)_\d(x_0) - (R(s)f)_\d(x_0)\big|\ge 2-\e,
$$
where $f_\d$ is defined as in the previous lemma. Since $f_\d\in
BUC^\circ(E)$, $\n f_{\d}\n\le 1$, and $\n R^\circ(t)\n\le 1$,
$\n R^\circ(s)\n \le  1$, the lemma follows. \epf

In combination with the technique described in Introduction we
obtain a similar result for the Ornstein-Uhlenbeck semigroup:

\bpr\label{prop:norm-disc-P-Cb} For all $t,s\ge 0$ such that
$S(t)\not=S(s)$ we have $\n
P^\circ(t)-P^\circ(s)\n_{\calL(BUC^\circ(E))} =2.$ \epr
 \bpf
Define the invertible isometries $V_n: BUC (E)\to
BUC(E)$ by
$$ V_n f(x) = f(n^{-1}x), \qquad x\in E, \ \ f\in BUC (E).
$$
We will show that  $ L_{R^\circ}$ belongs to the limit class of
$L_{P^\circ}$. To this end,  for any $f \in BUC (E)$ and $x \in E
$, one has
$$ | V^{-1}_n P(t) V_n f (x) - R(t) f (x)| \le
  \int_E \Big| f (S(t)x +  n^{-1}y) - f (S(t) x ) \Big|\, d\mu_t(y)
  \le
\int_E  \omega_f (n^{-1}y)\, d\mu_t(y),
$$
where $\omega_f$  denotes the modulus of continuity of $f$.
Letting $n \to \infty $,  the last term tends to 0 by the
dominated convergence theorem. Hence, for any $f \in BUC (E)$,
$$
\limn \| V^{-1}_n P(t) V_n f  - R(t) f \| =0.
$$
 The result now follows from Lemma
\ref{lem:norm-disc-R-Cb}.
\epf

\begin{corollary} If $A\not=0$, then there exists $t_0>0$ such that
$\n P^\circ(t)- P^\circ(s) \n_{\calL(BUC^\circ(E))} =2$ for all
$0\le t,s\le t_0$, $t\not=s$.
\end{corollary}
\begin{proof}
If such $t_0$ does not exist, there exist sequences $s_n\da 0$ and $t_n\da
0$ with $s_n\le t_n$ such that $\n P^\circ(t_n)- P^\circ(s_n) \n_{\calL(BUC^\circ(E))}<2$ for all
$n$. By Proposition \ref{prop:norm-disc-P-Cb}, 
$S(s_n) = S(t_n)$ for all $n$.  Fixing an element $x\in D(A)$, for all $n$ we obtain
$$ \int_{s_n}^{t_n} S(r)Ax\,dr = S(t_n)x-S(s_n)x = 0.$$ 
Upon dividing both sides by $t_n-s_n$ and passing to the limit $n\to\infty$
we obtain
$Ax = 0$. This being true for all $x\in D(A)$ we conclude that $A=0$.
\end{proof}

By a result of \cite{NZ} the same conclusion holds for $A=0$
if the range of $Q$ is infinite-dimensional; see also \cite{DR}
where the special case of a
Hilbert space $E$ was considered.

We proceed with a different sufficient condition for norm
discontinuity which, for the case of a Hilbert spaces $E$, is
implicitly contained in \cite{NZ}. Two probability measures $\mu$
and $\nu$ are called {\em disjoint}, notation $\mu\perp\nu$, if
there exist disjoint measurable sets $A$ and $B$ such that $\mu(A)
= \nu(B)=1$. The measures $\mu$ and $\nu$  are called {\em
equivalent}, notation $\mu\sim\nu$, if they are mutually
absolutely continuous, i.e., $\mu\ll\nu$ and $\nu\ll \mu$.

\bpr For all $t,s\ge 0$ such that $\mu_t\perp\mu_s$ we have $\n
P^\circ(t) - P^\circ(s)\n_{\calL(BUC^\circ(E))} = 2$. \epr \bpf By
assumption we have $\n \mu_t - \mu_s\n_{{\rm var}}=2$, where
$\n\cdot\n_{\rm var}$ denotes the total variation norm of a finite
signed Radon measure. Identifying $\mu_t$ and $\mu_s$ with
elements from the dual of $BUC^\circ(E)$, it will be enough to
show that $\n \mu_t - \mu_s\n_{(BUC^\circ(E))\s} = 2$. Indeed,
once we know this, given $\e>0$ we choose $g\in BUC^\circ(E)$ with
$\n g\n=1$ such that $|\lb g,\mu_t - \mu_s\rb|\ge 2-\e$ and
observe that
$$
 \n P^\circ(t) - P^\circ(s)\n_{\calL( BUC^\circ(E))} \ge
 |P^\circ(t)g(0) - P^\circ(s)g(0)| = |\lb g,\mu_t - \mu_s\rb|\ge 2-\e.
$$
 Suppose $\nu$ is a finite signed  Radon measure on $E$.
 Generalizing \cite[Lemma 2.3]{NZ}, the proof will be finished
by showing that \beq\label{eq:var} \n \nu\n_{(BUC^\circ(E))\s} =
\n \nu\n_{\rm var}. \eeq
The inequality `$\le$' is clear.  To
check the inequality `$\ge$', by the Jordan-Hahn decomposition it
is enough to prove the assertion when $\nu$ is a Radon probability
measure on $E.$ By \cite[Section 1.1]{Bi}, for any given $\e >0$
there exists $f\in BUC(E)$ with $\n f\n\le 1$ such that $\lb
f,\nu\rb\ge 1-\e$. For $\d>0$ define $f_\d\in BUC^\circ(E)$ as in
Lemma \ref{lem:delta}. By inner regularity of $\nu$, the supremum
of $\nu (K)$ with $K$ ranging over all compact subsets of $E$
equals $1$. Hence to prove \eqref{eq:var} it is enough to observe
that by Lemma \ref{lem:delta} we have $\lim_{\d\da 0}f_\d=f$
uniformly on compact sets. \epf  In the converse direction we have
the following result.

\bpr If $t,s\ge 0$ are such that $S(t)=S(s)$ and $\n P^\circ(t) -
P^\circ(s)\n_{\calL(BUC^\circ(E))}$
$  = 2$, then
$\mu_t\perp\mu_s$. \epr \bpf Given $\e>0$, there exist $f\in
BUC^\circ(E)$ and $x\in E$ such that $ |P^\circ(t)f(x) -
P^\circ(s)f(x)| \ge 2-\e.$ Defining $g\in BUC^\circ(E)$ by $g(y) =
f(S(s)x+y)$, this may be restated as
$$
\bal
\ &
\Big| \int_E g(y)\,d\mu_t(y) - \int_E g(y)\,d\mu_s(y) \Big|
\\ & \qquad\qquad  = \Big| \int_E f(S(s)x+y)\,d\mu_t(y) -\int_E f(S(s)x+y)\,d\mu_s(y) \Big|
\\ & \qquad\qquad  = \Big| \int_E f(S(t)x+y)\,d\mu_t(y) - \int_E f(S(s)x+y)\,d\mu_s(y) \Big|
\\ & \qquad\qquad = |P^\circ(t)f(x) - P^\circ(s)f(x)| \ge 2-\e.
\eal
$$
This shows that $\n \mu_t - \mu_s\n_{(BUC^\circ(E))\s} \ge 2-\e$.
Since the choice of $\e>0$ is arbitrary we obtain that
$$ 2\le \n \mu_t - \mu_s\n_{(BUC^\circ(E))\s}\le
\n \mu_t - \mu_s\n_{{\rm var}} \le  2,$$
the second and third of these inequalities being obvious.
Hence $\n \mu_t - \mu_s\n_{{\rm var}}=2$,
which implies that $\mu_t\perp\mu_s$.
\epf

By putting these results together we have proved:

\bth\label{thm:mainBUC} For all $t,s\ge 0$ the following
assertions are equivalent: \ben \item $\n P^\circ(t) -
P^\circ(s)\n_{\calL(BUC^\circ(E))}
= 2$; \item $S(t)\not=S(s)$ or
$\mu_t\perp\mu_s$. \een \eth

It should be observed that neither  $S(t)\not=S(s)$ implies
$\mu_t\perp\mu_s$ or conversely. An example of a periodic
semigroup with period $1$ such that $\mu_t\perp\mu_s$ for all
$t,s\ge 1$ is given in \cite[Example 3.8]{NeJFA}. On the other
hand, if $\dim E<\infty$, then for any choice of $S$ and $B$ the
measures $\mu_t$ and $\mu_s$ are mutually absolutely continuous
for all $t,s\ge t_0$.

We continue with two examples which show that $\n P(t) -
P(s)\n_{\calL(BUC^\circ(E))} < 2$ may occur for  certain values of
$t\not=s$.

\bex[Nilpotent $S$] Let $E=L^2(0,1)$ and let $S$ be the nilpotent
shift semigroup on $L^2(0,1)$, see for instance \cite[page 120]{EN}.
Then $S(t) = S(s)=0$ and $\mu_t =
\mu_s = \mu_1$  for all $t,s\ge 1$. Hence, $P(t)=P(s)$ for all
$t\ge s\ge 1$. \eex

\bex[Periodic $S$ in finite dimensions] Let $ H=E=\R^2$ and let $S$
be the rotation group on $\R^2$. Let $B:= I$, the identity
operator on $\R^2$.  Since $S\s(t) = S(-t)$ for all $t\ge 0$, the
covariance operator of $\mu_t$ is given by $Q_{t} = t I.$ Hence
$\mu_t$ is the Gaussian measure on $\R^2$ with variance $t$. For
$k=0,1,2,\dots$ and $f\in BUC^\circ(\R^2)$,
$$ P^\circ(2k\pi)f(x) = \int_E f(S(2k\pi)x+y)\,d\mu_{2k\pi}(y)
=  \int_E f(x+y)\,d\mu_{2k\pi}(y).
$$
For $j\ge 1$, $k\ge 1$, $j\not=k$, we have $S(2j\pi)=S(2k\pi)$ and
$\mu_{2j\pi}\sim \mu_{2k\pi}$. Theorem \ref{thm:mainBUC}
shows that $\n P^\circ(2j\pi) - P^\circ(2k\pi)\n < 2.$ \eex

We will show next that the above two examples are in some
sense the only possible ones.

Recall that a Gaussian measure $\nu$ on $E$ is called {\em nondegenerate}
if there exists no proper closed subspace $E_0$ of $E$ with $\nu(E_0)=1$. It is easy to see that $\nu$ is nondegenerate if and only if its covariance operator has dense range.

For $t>0$ fixed, $P$ is said to be {\em
strongly Feller at time $t$} if $P(t)f \in C_b(E)$ for all $f\in
B_b(E)$. Here $B_b(E)$ denotes the space of real-valued bounded
Borel functions on $E$. As is well known, $P$ is strongly Feller
at time $t$ if and only if we have $S(t)E \subseteq H_{Q_t}$,
where $H_{Q_t}$ is the reproducing kernel Hilbert space associated
with $Q_t$;  cf. \cite{DZ, NeJFA}.

\bth \label {ciao} Suppose $t>s\ge 0$  are such that
$\n P^\circ(t) - P^\circ(s)\n_{\calL(BUC^\circ(E))} < 2.$
Assume in addition that one of the following two assumptions is satisfied:
\begin{enumerate}
\item[(i)] $\mu_{t-s}$ is nondegenerate;
\item[(ii)] $P$ is strongly Feller at time $t-s$.
\end{enumerate}
Then there exists a direct sum decomposition into $S$-invariant
subspaces $E = E_0\oplus E_1$, with $\dim E_1 <\infty$, such that
$S$ is nilpotent on $E_0$ and periodic on $E_1$ with period $t-s$.
\eth

\bpf
By Theorem \ref{thm:mainBUC}, the assumption of the theorem implies that
$S(t)=S(s)$ and $\mu_t\not\perp \mu_s$. By the Feldman-Hajek theorem
\cite[Theorem 2.7.2]{Bo},
$\mu_t\sim \mu_s$.

Let $H_Q$ be the reproducing kernel Hilbert space associated
with $Q =  B B^*$ and let $E_s$ denote the closure of the range of
$S(s)$.  Define $j:H_Q\to E_s$ by $j := S(s)B$ and
$R\in\calL(E_s\s,E_s)$ by $R:= jj\s = S_s(s)QS_s^*(s)$, where
$S_s(s)$ is the operator $S(s)$ regarded as an operator from $E$
to $E_s$.
 For $\tau>0$ introduce  the operators $R_\tau\in\calL(E_s\s,E_s)$ by
$$
R_\tau y\s := \int_0^\tau S(u)RS\s(u)y\s\,du, \qquad y\s\in E_s\s,$$
where by abuse of notation we think of $S$ as a semigroup on $E_s$.
Then $R_\tau$ is the covariance operator of the image measure
$\nu_\tau = S_s(s)\mu_\tau$ on $E_s$, i.e., $R_\tau  = S_s(s) Q_{\tau}
S_s\s(s).$ Moreover,
\beq\label{eq:equiv}
 \nu_{s} = S_s(s)\mu_{s} \sim S_s(s)\mu_{t} = \nu_{t}. \eeq Clearly,
\beq\label{eq:per} S(t-s)|_{E_s} = I|_{E_s}.
 \eeq
 By \eqref{eq:equiv} and \cite[Corollary 3.3]{NeJFA}, for
$\tilde k \in \N$ such that $\tilde k (t-s) \ge s$ we obtain
  $$\nu_{\tilde k  (t-s)} =\nu_{ s + (\tilde k  (t-s) - \, s )}  \sim
 \nu_{ t + (\tilde k  (t-s) - \, s )} = \nu_{(\tilde k +1)(t-s)}.$$
 But by \eqref{eq:per} we have $R_{(\tilde k + 1) (t-s)} = (\tilde k
+1)R_{t-s}$, and therefore the Feldman-Hajek theorem implies that
the reproducing kernel Hilbert space $H_{R_{t-s}}$ associated with
$R_{t-s}$ is finite-dimensional; cf. \cite[Example 2.7.4]{Bo}.

We will show below that each of the conditions (i) and (ii)
implies that the measure $\nu_{t-s}$ is nondegenerate.
Once we know this, the
proof can be finished as follows. Since $\nu_{t-s}$ is
nondegenerate, the reproducing kernel Hilbert space $H_{R_{t-s}}$
is dense in $E_{s}$. It follows that $H_{R_{t-s}} = E_s$, which
means that $E_s$ is finite-dimensional. Hence $E_s$ equals the
range of $S(s)=S(t)$. By the semigroup property, $E_s$ equals also
the range of $S(k(t-s))$, where the integer $k\ge 1$ is such that
$s\le k(t-s)< t$. In combination with \eqref{eq:per} it follows
that $S(k(t-s))$ is a projection in $E$. This proves the theorem,
with $E_0: = \ker S(k(t-s))$ and $E_1:=E_s = {\rm
range}\,S(k(t-s))$.

{\vskip 1mm}To finish the proof we show that both (i) and (ii)
imply the nondegeneracy of the measure $\nu_{t-s}$.

First assume (i).
It is immediate from the definition that the image of a nondegenerate
Gaussian measure under a bounded operator with dense range is nondegenerate.
Thus the nondegeneracy
assumption on $\mu_{t-s}$ implies that $\nu_{t-s}$ is
nondegenerate.

Next assume (ii).
Write $H_{t-s}:= H_{Q_{t-s}}$ for brevity and let $i_{t-s}:
H_{t-s}\embed E$ be the inclusion mapping. Recalling that $Q_{t-s}
= i_{t-s}\circ i_{t-s}\s$, for all $u\s\in E_s\s$ and $x\s\in E\s$
such that $x\s|_{E_s} = u\s$ we have
$$
 \langle R_{t-s} u\s,u\s \rangle
= \langle S_s(s)Q_{t-s} S_s\s(s) u\s,u\s \rangle = \langle
S(s)Q_{t-s} S\s(s) x\s, x\s \rangle = \n i_{t-s}\s S\s(s)x\s
\n_{H_{t-s}}^2.
$$
By the strong Feller property and a closed graph argument,
$S(t-s)$ is bounded as an operator from $E$ to $H_{t-s}$. Denoting
this operator by $T(t-s)$ we have $S(t-s) = i_{t-s}\circ T(t-s)$.
Let $I: E_s\to E$ be the inclusion mapping.
On $E_s$ we have $S(t)\circ I = I\circ S(t)$, where as before we
abuse of notation by writing $S$ for the restriction of $S$ to
$E_s$. Then, for all $x\s\in E\s$,
$$
\bal
 \n S\s(t)I\s x\s\n  & =\n I\s S\s(t-s)S\s(s) x\s\n
 \\ & = \n I\s T\s(t-s) i_{t-s}\s S\s(s)x\s\n
 \le \n T(t-s)I\n_{\calL(E_s, H_{t-s})}\,\n i_{t-s}\s S\s(s)x\s\n_{H_{t-s}}.
\eal
$$
Combining these things we obtain
$$\n T(t-s)I\n_{\calL(E_s, H_{t-s})}^{2}\, \langle R_{t-s} I\s x\s,I\s x\s \rangle \ge
\n S\s(t)I\s x\s\n^2 \ge c_{t}^2 \n I\s x\s\n^2, \qquad x\s\in
E\s,
$$
where the last estimate follows from the fact that $S$ is periodic
on $E_s$. Since $I\s$ is a surjection from $E\s$ onto $E_s\s$,
this gives that either  $R_{t-s}$ is nondegenerate or $T(t-s)I=0$.
In the first case the proof is complete.
If $T(t-s)I=0$, then $S(t-s)I=0$ as well, which
means that $S(t-s)=0$ on $E_s$. By periodicity this implies that
$E_s  = \{0\}$. This in turn implies that $S(s)=0$, i.e., $S$ is
nilpotent on $E$.
\epf

The nondegeneracy assumption on $\mu_{t-s}$  in (i) is fulfilled
if $Q$ has dense range; this is proved in the same way as
\cite[Lemma 5.2]{GGN}.

\begin{corollary} Let $\dim E= \infty$, and assume that $S$ is analytic
and condition {\rm (i)} or {\rm (ii)} is satisfied. Then
for all $t>s\ge 0$ we have $ \n P^\circ(t)-P^\circ(s)\n_{BUC^\circ(E)} = 2.$
\end{corollary}
\begin{proof}
An analytic $C_0$-semigroup on a nonzero Banach space cannot be nilpotent.
Hence, Theorem \ref{ciao}
shows that if there exist $t>s\ge 0$ such that $\n P^\circ(t)-P^\circ(s)\n_{BUC^\circ(E)} < 2$,
then $\dim E< \infty$.
\end{proof}

Next we consider the case $A=0$. In this situation one has $Q_t =
tQ = tBB\s$, and since by our standing assumption these operators
are Gaussian covariances, it follows that $Q$ is a Gaussian
covariance. We denote the Gaussian measure on $E$ with covariance
operator $Q$ by $\nu$. The semigroup $P$ is then the heat
 semigroup
  given by
$$P(t)f(x) =\int_E f(x+y)\,d\mu_t(y) = \int_E f(x+\sqrt{t}\,y)\,d\nu(y),
\qquad f\in C_b(E).
$$ The
 restriction of $P$ to $BUC(E)$ is
strongly continuous with respect to the supremum norm. The
infinitesimal generator $L_P$ of $P$ is given, on a suitable core
of cylindrical functions, by
$$L_{P} f(x) = \tfrac12\hbox{Tr}\,QD^2f(x).$$
The following result was proved in \cite{MRS}  for the
special case of an infinite-dimensional Hilbert space $E$.
Our proof is essentially the same, the main difference being that
the coordinate-free presentation simplifies matters somewhat. The
spectrum and approximate point spectrum of $L_P$ in $BUC(E)$ are
denoted by $\sigma(L_P)$ and $\sigma_{\rm a}(L_P)$, respectively.

\begin{proposition}
If $A=0$ and $Q$ is not of finite rank, then
$\sigma(L_P) = \sigma_{\rm a}(L_P) = \C^-.$
\end{proposition}
\begin{proof}
Fix a sequence $(x_n\s)$ in $E\s$ such that $(B\s x_n\s)$ is an
orthonormal sequence in $H$. Such a sequence exists since
the range of $B\s$ is not finite-dimensional in $H$.
For each $n\ge 1$ we consider the map
$T_n: E\to \R^n$ defined by
$$T_nx := (\langle x,x_1\s\rangle,\dots , \langle x,x_n\s\rangle).$$
The image measure of $\nu$ under $T_n$ equals $\gamma_n$, the standard
Gaussian measure on $\R^n$.

Let $\Delta_n$ be the Laplace operator acting in $BUC(\R^n)$.
Denoting the heat semigroup on $BUC(\R^n)$
generated by $\tfrac12\Delta_n$ by $\{P_n(t)\}_{t\ge 0}$,
for all $f\in BUC(\R^n)$ and $x\in E$ we have
$$
P(t) f(T_n x)
 = \int_E f(T_n(x+ \sqrt{t}\,y))\,d\nu(y)
 = \int_{\R^n} f(T_n x+\sqrt{t}\,\eta ))\,d\gamma_n(\eta)
 = P_n(t)f(T_n(x)).
$$
From this it is immediate that $f\circ T_n \in D(L_P)$ whenever $f\in
D(\Delta_n)$
and  in this case,
$$L_P (f\circ T_n) = (\tfrac12\Delta_nf)\circ T_n.$$
Fix $\l \in\C$ with $\Re\l < 0$ and
consider the functions $f_{n,\l}, g_{n,\l}\in BUC(\R^n)$ defined
by
$$
f_{n,\l}(\xi) = \exp\Bigl(\frac{\l}{n}|\xi|^2\Bigr)\quad
\hbox{and}\quad g_{n,\l}(\xi) =
\frac{- 2 \lambda^2 |\xi|^2}{n^2}f_{n,\l}(\xi),\qquad \xi\in\R^n.
$$
We have $f_{n,\l}\in
D(\Delta_n)$ and
$$(\l-\tfrac12\Delta_n) f_{n,\l}= g_{n,\l}.$$
Hence $f_{n,\l}\circ T_n\in
D(L_P)$ and
$$(\l-L_P)(f_{n,\l}\circ T_n) = g_{n,\l}\circ T_n.$$
Moreover,
$$\Vert f_{n,\l}\circ T_n\Vert_{BUC(E)} = \Vert f_{n,\l}\Vert_{BUC(\R^n)} = 1$$
and we compute
$$\Vert g_{n,\l}\circ T_n\Vert_{BUC(E)} = \Vert g_{n,\l}\Vert_{BUC(\R^n)}
= \frac{2|\l|^2}{ne\,|\Re\l|}.
$$
This proves that the sequence
$(f_{n,\l}\circ T_n)$ is an approximate eigenvector for $L_P$,
with approximate eigenvalue $\l$.
It follows that $\{\Re\l<0\}\subseteq \sigma_{\rm a}(L_P)$. On the other hand,
since $\{P(t)\}_{t\ge 0}$ is a contraction semigroup on $BUC(E)$,
we have $\{\Re\l>0\}\subseteq\varrho(L_P)$, where $\varrho(L_P)$
denotes the resolvent set of $L_P$. Combining this, we see that
$\sigma(L_P) = \C^-$. Moreover,
$i\R = \partial\sigma(L_P)\subseteq \sigma_{\rm a}(L_P)$ by the general theory of semigroups, and therefore
$\sigma(L_P) = \sigma_{\rm a}(L_P) = \C^-$.
\end{proof}

If $A=0$ and $E=\R^d$, then $P$ is analytic and therefore
$\sigma(L_P)$ is contained in a strict subsector in $\C^-$. For
$A\not=0$, $Q $ invertible,  and $E=\R^d$, it was shown in
\cite{M} that $\sigma(L_{P^\circ}) =\C^-$ if $\sigma(A)\subseteq
\{\l\in\C: \ \Re\l<0\}$ and that $\sigma(L_{P^\circ}) \supseteq 
\C^-$ if $\sigma(A)\subseteq
\{\l\in\C: \ \Re\l>0\}$, and that in both cases 
every $\l\in\C$ with $\Re\l<0$ is an
eigenvalue. We have the following extension of this result
to infinite dimensions:

\begin{theorem} Assume that the operator $Q$ has  dense range.
Assume also that $S$ is eventually compact and that $\sigma (A)$
is nonempty.

\begin{enumerate}
\item If $\sigma(A)\subseteq \{\l\in\C: \ \Re\l<0\}$, 
then $\sigma(L_{P^\circ}) =\C^-$ and every $\l\in\C$ with $\Re\l<0$ 
is an eigenvalue.
\item If $\sigma(A)\subseteq \{\l\in\C: \ \Re\l>0\}$, 
then $ \sigma(L_{P^\circ})\supseteq\C^-$ and every $\l\in\C$ with $\Re\l<0$ 
is an eigenvalue.
\end{enumerate}
\end{theorem}

The proof is based on the same Riesz projection
argument as Theorem \ref{cor:sp-L1} below  and is left to the reader.

\section{The Ornstein-Uhlenbeck semigroup in spaces of integrable functions}
\label{sec:L1}

Our approach to proving norm discontinuity of Ornstein-Uhlenbeck
semigroups in $L^1$-spaces is based on the following observation.

\ble\label{lem:disj}
For all $t>s\ge 0$ with $S(t-s)\not=I$ there exist $x_0\in E$ and $r>0$ such
that
$$ \{x\in E: \ \n S(t)x - x_0 \n < r\} \cap  \{x\in E: \ \n S(s)x - x_0 \n <
r\} = \emptyset.$$
\ele
\bpf
Choose $x_0\in E$ such that $S(t-s)x_0\not=x_0$.
Let $M := \max\{1,\n S(t-s)\n_{\calL(E)}\}$ and put
$$ r:= \frac1{2M}\n S(t-s)x_0-x_0\n.$$
Suppose $x\in E$ is such that $\n S(s)x - x_0 \n < r$.
We will prove that $\n S(t)x - x_0 \n \ge  r$.
By assumption there exists a vector $x_1\in E$ with $\n x_1\n< r$ such that
$S(s)x = x_0+x_1$. Then,
$$
\bal
\n S(t)x-x_0\n & = \n S(t-s)(x_0+x_1) -x_0 \n
\\ & \ge \n S(t-s)x_0-x_0\n - \n S(t-s)x_1\n \ge 2Mr - Mr  = Mr \ge r.
\eal$$
\epf

Until further notice we now specialize to the case where
$E=\R^d$ and assume that $A$ is an $(d\times d)$-matrix with real
coefficients. We write $S(t) = e^{tA}$. As before, $R$
indicates the drift semigroup given by \eqref{3}. Let $C_c(\R^d)$
denote the space of continuous compactly supported functions on
$\R^d$.

For all $f\in C_c(\R^d)$ we have \beq\label{eq:det}
\int_{\R^d} |R(t)f(x)|\,dx
= \frac{1}{|\det(S(t))|} \int_{\R^d} |f(S(t)x)|\, |\det(S(t))|\,dx
= e^{-t{\rm Tr}\,A} \int_{\R^d} |f(y)|\,dy.
\eeq
It follows that the
restrictions of $R(t)$ to $C_c(\R^d)$ extend to bounded operators
on $L^1(\R^d)$ of norm $ \n R(t)\n_{L^1(\R^d)} = e^{-t{\rm
Tr}\,A}.$ Since also $\lim_{t\downarrow 0}\n R(t)f-f\n_{L^1(\R^d)}
= 0$ for all $f\in C_c(\R^d)$ it follows that $R$ has a unique
extension to a $C_0$-semigroup on $L^1(\R^d)$. The space
$C_c^1(\R^d)$ is a core for its generator $L_R$
 and we have
$$
L_R f(x) = \lb Ax, D f(x)\rb, \qquad x\in\R^d, \
f\in C_c^1(\R^d).
$$

\bpr\label{prop:norm-disc-R} For all $t>s\ge 0$ with
$S(t)\not=S(s)$ we have $ \big\n e^{t{\rm Tr}\,A}R(t)-e^{s{\rm
Tr}\,A}R(s)\big\n_{\calL(L^1(\R^d))} =2.$
\epr \bpf
Let
$x_0\in\R^d$ and $r>0$ be as in Lemma \ref{lem:disj}. By Lemma
\ref{lem:disj},
$$
\bal
\ &  \big\n e^{t{\rm Tr}\,A}R(t)\one_{\{\n x-x_0\n<r\}}-
e^{s{\rm Tr}\,A}R(s)\one_{\{\n x-x_0\n<r\}}\big\n
\\ & \quad\quad = \big\n  e^{t{\rm Tr}\,A}\one_{\{\n S(t)x-x_0\n<r\}}
- e^{s{\rm Tr}\,A}\one_{\{\n S(s)x-x_0\n<r\}}\big\n
\\ & \quad\quad = e^{t{\rm Tr}\,A}\n  \one_{\{\n S(t)x-x_0\n<r\}}\n
+ e^{s{\rm Tr}\,A}\n\one_{\{\n S(s)x-x_0\n<r\}}\n
\\ & \quad\quad =e^{t{\rm Tr}\,A}\n R(t) \one_{\{\n x-x_0\n<r\}}\n
+ e^{s{\rm Tr}\,A}\n R(s)\one_{\{\n x-x_0\n<r\}}\n
= 2\n \one_{\{\n x-x_0\n<r\}}\n,
\eal$$
where in the last step we used \eqref{eq:det}.
It follows that
$ \big\n e^{t{\rm Tr}\,A}R(t)-e^{s{\rm Tr}\,A}R(s)\big\n \ge 2.$
Since by \eqref{eq:det} we also have
$e^{\tau{\rm Tr}\,A}\n R(\tau)\n\le 1$ for all $\tau\ge 0$, the proposition
is proved.
\epf

Our next aim is to extend the Ornstein-Uhlenbeck semigroup $P$ to
$L^1(\R^d)$ as well. For all $f\in C_c(\R^d)$ we have
\beq\label{eq:detP} \bal \int_{\R^d} |P(t)f(x)|\,dx & \le
\int_{\R^d}\int_{\R^d} |f(S(t)x +y)|\,dx\,d\mu_t(y)
\\ & = e^{-t{\rm Tr}A} \int_{\R^d}\int_{\R^d} |f(\xi)|\,d\xi\,d\mu_t(y)
 = e^{-t{\rm Tr}A} \int_{\R^d}|f(\xi)|\,d\xi
\eal \eeq
with equality for nonnegative functiones $f$.
 It follows that the restrictions of the operators $P(t)$
to $C_c(\R^d)$ extend to bounded operators on $L^1(\R^d)$ of norm
$ \n P(t)\n_{\calL( L^1(\R^d))} =e^{-t{\rm Tr}A}.$ Since also
$\lim_{t\downarrow 0}\n P(t)f-f\n_{L^1(\R^d)} = 0$ for all $f\in
C_c(\R^d)$ it follows that the restriction of $P$ to $C_c(\R^d)$
has a unique extension to a $C_0$-semigroup on $L^1(\R^d)$, which
is still given by formula \eqref{v}. The space $C_c^2(\R^d)$ is a
core for its generator $L_P$ and we have
$$
L f(x) =\tfrac{1}{2} {\rm Tr}\, QD^2 f(x) + \lb Ax, Df(x)\rb, \qquad x\in\R^d, \
f\in C_c^2(\R^d).
$$

\bth\label{thm:norm-disc-P} For all $t>s\ge 0$ with
$S(t)\not=S(s)$ we have $\big\n e^{-t{\rm Tr}\,A}P(t)-e^{-s{\rm
Tr}\,A}P(s)\big\n_{\calL(L^1(\R^d))} =2.$ \eth \bpf For
$n=1,2,\dots$ let $V_n: L^1(\R^d)\to L^1(\R^d)$ denote the
invertible isometry
$$ V_n f(x) = n^{-d} f(n^{-1}x), \qquad x\in \R^d, \ f\in L^1(\R^d).$$
As in the proof of Proposition \ref{prop:norm-disc-P-Cb} we see that
$L_R$ belongs to the limit class of $L_P$.
Hence by Proposition \ref{prop:norm-disc-R} and \cite[Proposition 12]{DS}, applied to the operators $L_P-{\rm Tr}\,A$ and
$L_R-{\rm Tr}\,A$,
$$  \big\n e^{-t{\rm Tr}\,A}P(t)-e^{-s{\rm Tr}\,A}P(s)\big\n_{\calL(L^1(\R^d))}
\ge  \big\n e^{-t{\rm Tr}\,A}R(t)-e^{-s{\rm Tr}\,A}R(s)\big\n_{\calL(L^1(\R^d))}
\ge 2.
$$
Since by \eqref{eq:detP} we also have
$e^{\tau{\rm Tr}\,A}\n P(\tau)\n\le 1$ for all $\tau\ge 0$, the theorem
is proved.
\epf

Alternatively this theorem may be derived from Proposition 
\ref{prop:norm-disc-P-Cb} via the
duality argument of \cite[Lemma 3.6]{M}. 


After these preparations  we come to the main results of this
section, which give conditions for norm discontinuity of $P$ in
the space $L^1(E,\mu_\infty)$, where $\mu_\infty$ is the invariant
measure for $P$ discussed in Section \ref{sec:intr}. Note that in
finite dimensions, its existence is guaranteed under the mere
assumption that the limit $Q_\infty:=\lim_{t\to\infty} Q_t$ exists
in $\calL(\R^d)$. This will be assumed in the next result,
in which $P$ denotes Ornstein-Uhlenbeck semigroup on
$L^1(\R^d,\mu_\infty)$ and $L_P$ its generator. Since we are
dealing with the finite-dimensional case, a sufficient condition
for the existence of $Q_\infty$ is that $\sigma(A)\subseteq
\{\Re\l<0\}$.

\bth\label{thm:main} Assume that the limit
$Q_\infty:=\lim_{t\to\infty} Q_t$ exists in $\calL(\R^d)$ and let
$\mu_\infty$ be the Gaussian measure on $\R^d$ with covariance
matrix $Q_\infty$. Then for all $t,s\ge 0$ with $t\not=s$ we have
$$  \n P(t)-P(s)\n_{\calL(L^1(\R^d,\mu_\infty))} =2.$$
\eth
\bpf
As is well known \cite[Proposition 1]{CG}, the range of $Q_\infty$
is invariant under the action of $S$ and therefore
we may assume without loss of generality that $\mu_\infty$ is nondegenerate.
Moreover, the existence of $\mu_\infty$ implies that $S(t)\not=S(s)$ for all
$t,s\ge 0$ with $t\not=s$, since otherwise the improper integral
defining $Q_\infty$ will diverge.

Let $b$ be the density of $\mu_\infty$ with respect to
the Lebesgue measure; this
density exists since $\mu_\infty$ is assumed to be nondegenerate.
Proceeding as in \cite{MPP} we consider the invertible isometry
$V: L^1(\R^d) \to L^1(\R^d,\mu_\infty)$ given by $f\mapsto
b^{-1}f$ and define $\tilde P(t): L^1(\R^d)\to
L^1(\R^d)$ by
$$ \tilde P(t) = V^{-1}\circ P(t)\circ V, \qquad t\ge 0.$$
Then $\tilde P = \{\tilde P(t)\}_{t\ge 0}$ is a $C_0$-semigroup on
$L^1(\R^d)$ and by the computations in \cite[Theorem 5.1]{MPP} its
generator $\tilde L$ is given by
 \beq\label{eq:tildeL} \tilde L
f(x)= \tfrac12 {\rm Tr}\, QD^2 f (x) + \lb \tilde Ax,D f(x)\rb +
kf(x),\qquad x\in \R^d, \ \ f\in C_c^2(\R^d),
 \eeq where
$$
\tilde A= -Q_\infty A\s Q_\infty^{-1}, \quad k
= -{\rm Tr}\,A = -{\rm Tr}\,\tilde A.
$$
The result now follows
from Theorem \ref{thm:norm-disc-P} applied to $\tilde L - k$. \epf
   
Returning to the setting of an arbitrary real Banach space $E$,
we have the following extension of Theorem \ref{thm:main}.

\begin{theorem}\label{thm:ev-comp}
Assume that the weak operator limit $Q_\infty:=\lim_{t\to\infty}
Q_t$ exists in $\calL(E\s,E)$ and that it is the covariance
operator of a Gaussian measure $\mu_\infty$ on $E$.  Let $S$ be an
eventually compact $C_0$-semigroup on $E$, and assume that its
generator $A$ has nonempty spectrum.
Then for all $t,s\ge 0$ with $t\not=s$ we have
 \beq \label{f2}
  \n P(t)-P(s)\n_{\calL(L^1(E,\mu_\infty))} =2.
 \eeq
\end{theorem}
\bpf Replacing $E$ by the closure of the reproducing kernel space
associated with $Q_\infty$, which is invariant under $S$ by
\cite[Proposition 1]{CG}, see also \cite{NeJFA}, we may assume
without loss in generality that $\mu_\infty$ is nondegenerate.

Since $\sigma(A)\not=\emptyset$ we may fix some $\lambda_0
\in \sigma(A)$. Note that $\lambda_0$ is an isolated point in
$\sigma (A)$.  Let $\pi_0:E\to E$ be the Riesz
projection onto $E_0$, the finite dimensional subspace of $E$
generated by all generalized eigenvectors associated to
$\lambda_0$, cf. \cite[Corollary 3.2, page 330]{EN}. The
projection $\pi_0$ commutes with the operators $S(t)$. Let $S_0$
denote the restriction of $S$ to $E_0$, with generator $A_0
\in \calL(E_0)$, and define $Q_0\in\calL(E_0\s,E_0)$ by $Q_0 :=
\pi_0 Q \pi_0\s$. Here we think of $\pi_0$ as an operator from $E$ onto
$E_0$. For $0\le t\le \infty$ the covariance operator $Q_{0,t}$
associated with the image measure $\mu_{0,t} = \pi_0 \mu_t$ on
$E_0$ is given by
$$
 Q_{0,t} x_0\s = \int_0^t S_0(s) Q_0  S_0\s(s) x_0\s \, ds,\qquad
x_0\s \in E_0\s.
$$
Since   $Q_{0, \infty}$ is nondegenerate and
$\sigma (A_0) = \{\lambda_0\}$ we have $\Re\lambda_0<0$.

For all $f\in L^1 (E_0, \mu_{0,\infty})$, the function
$$f_0(x):= f(\pi_0 x), \qquad x\in E, $$ belongs to
$L^1(E,\mu_{\infty})$ and we have \beq\label{eq:equal-norms}
 \int_E |f_0 (x )|\, d \mu_{\infty} (x) =
 \int_{E_0} |f (\xi)| \,d\mu_{0,\infty} (\xi).
\eeq
Let $P_0$ be the corresponding Ornstein-Uhlenbeck semigroup on $E_0$,
i.e.,
$$ P_0(t)f(x_0) = \int_{E_0} f(S_0(t)x_0+\xi)\,d\mu_{0,t}(\xi),
 \qquad t\ge 0, \ \ x_0\in E_0, \ \ f\in L^1(E_0,\mu_{0,\infty}).$$
With these notations,
\beq\label{eq:0} (P_0(t)f) (\pi_0 x) = P(t)f_0(x).\eeq
Now let $t>s\ge 0$ be such that \eqref{f2} holds. Then,
by virtue of \eqref{eq:equal-norms} and \eqref{eq:0},
$$
\bal \| P(t)  - P(s) \|_{\calL(L^1(E, \mu_{\infty}))} & \ge
\sup_{\|f\|_{L^1 (E_0, \mu_{0,\infty}) } \le 1 } \| P(t) f_0 -
P(s) f_0 \|_{L^1(E, \mu_{\infty})}
\\ & = \sup_{\|f \|_{L^1 (E_0, \mu_{0,\infty}) } \le 1 }
\| P_0(t) f - P_0(s) f \|_{L^1(E_0,  \mu_{0,\infty})}=2,
 \eal
$$
where the last step follows from the previous theorem. Since $P$
is a contraction semigroup in $L^1 (E,\mu_{\infty})$, the equality
\eqref{f2} follows.
 \epf

Our final result concerns the spectrum of $L_P$. The following description
of $\sigma(L_P)$ in $L^1(\R^d,\mu_\infty)$ was shown in
\cite{MPP}, where it was derived from the characterization of
$\sigma(L_P)$ for $L^1(\R^d)$, see  \cite{M}.

\begin{theorem}\label{thm:sp-L1} Assume that the limit
$Q_\infty:=\lim_{t\to\infty} Q_t$ exists in $\calL(\R^d)$ and let
$\mu_\infty$ be the Gaussian measure on $\R^d$ with covariance
matrix $Q_\infty$.
The spectrum of $L_P$ in $L^1(\R^d,\mu_\infty)$ equals
$\C^-$, and every $\l\in\C$ with $\Re\l < 0$ is an
eigenvalue of $L_P$.
\end{theorem}
In setting of a real Banach space $E$ we obtain the following extension:

\begin{theorem}\label{cor:sp-L1} Assume that the weak operator limit
$Q_\infty:=\lim_{t\to\infty} Q_t$ exists in $\calL(E\s,E)$ and that it is the
covariance operator of a Gaussian measure $\mu_\infty$ on $E$.
If $S$ is eventually compact and $\sigma(A)\not=\emptyset$,
the spectrum of $L_P$
in $L^1(E,\mu_\infty)$ equals $\C^-$,
and every $\l\in\C$ with $\Re\l < 0$ is an eigenvalue of $L_P$.
\end{theorem}
\bpf We may assume that $\mu_\infty$ is
nondegenerate. Fix $\l_0\in \sigma(A)$. Using the notations of the proof of Theorem \ref{thm:ev-comp}, let $L_{P_0}$ denote the generator of the
semigroup $P_0$ on $L^1(E_0,\mu_{0,\infty})$. Theorem \ref{thm:sp-L1} implies that $\sigma(L_{P_0})  = \C^-$ and that every $\l\in\C$ with $\Re\l <
0$ is an eigenvalue of $L_{P_0}$.
Let $f\in L^1(E_0,\mu_{0,\infty})$ be an associated eigenvector. Then
$f_0(x):= f(\pi_0 x)$  defines a function $f\in L^1(E,\mu_\infty)$ satisfying
$$P(t)f_0(x) = P_0(t)f(\pi_0 x) = e^{\l t} f(\pi_0 x) = e^{\l t} f_0(x).$$
Hence, $P(t)f_0 = e^{\l t} f_0$, and $f_0$ is an eigenvector for $L_P$ with
eigenvalue $\l$.
\epf

After the completion of this paper, the authors received the preprint \cite{Ch} by
Chojnowska-Michalik. She proves a related extension of
Theorem \ref{thm:sp-L1}: if the part of $A$ in the reproducing kernel Hilbert
space of $\mu_\infty$ has an
eigenvalue $\l$ with $\Re\l<0$,
then  $\sigma(L_P) = \C^-$.

\medskip \noindent{\em Acknowledgment.} \  The second author
wishes to thank the Delft Institute of Applied Mathematics,
 where this  paper was started, for hospitality and nice working
 conditions.

\end{document}